\documentclass{amsart}
\usepackage[english]{babel}
\usepackage[cp1250]{inputenc}
\usepackage{amsmath,amssymb}
\usepackage{microtype}
\usepackage{float}
\usepackage{enumitem}
\usepackage{mathrsfs}

\usepackage{tikz}
\usetikzlibrary{matrix,arrows}
\usetikzlibrary{decorations.markings}

\newcommand\blfootnote[1]{%
  \begingroup
  \renewcommand\thefootnote{}\footnote{#1}%
  \addtocounter{footnote}{-1}%
  \endgroup
}

\newtheorem{thm}{Theorem}[section]
\newtheorem{prop}[thm]{Proposition}
\newtheorem{lem}[thm]{Lemma}
\newtheorem{cor}[thm]{Corollary}

\newtheorem*{claim*}{Claim}

\theoremstyle{definition}

\newtheorem{remark}[thm]{Remark}
\newtheorem{example}[thm]{Example}

\newtheorem*{remark*}{Remark}

\def\FF{\mathbb F}

\def\NN{\mathbb N}
\def\ZZ{\mathbb Z}

\DeclareMathOperator{\tb}{tb}

\DeclareMathOperator{\rot}{rot}

\DeclareMathOperator{\HFK}{HFK}
\DeclareMathOperator{\CFK}{CFK}
\DeclareMathOperator{\gCFK}{gCFK}
\DeclareMathOperator{\HFKhat}{\widehat{HFK}}
\DeclareMathOperator{\SFH}{SFH}
\DeclareMathOperator{\HFhat}{\widehat{HF}}

\DeclareMathOperator{\Leg}{\mathfrak{L}}
\DeclareMathOperator{\Leghat}{\widehat{\mathfrak{L}}}
\DeclareMathOperator{\EH}{EH}
\DeclareMathOperator{\EHdirect}{\underrightarrow{\EH}}
\DeclareMathOperator{\EHinverse}{\underleftarrow {\EH}}

\begin{document}

\title[On contact surgery and knot Floer invariants]{On contact surgery and knot Floer invariants}

\author{Irena Matkovi\v{c}}
\address{Department of Mathematics, Uppsala University, Sweden}
\email{irena.matkovic@math.uu.se}

\begin{abstract}
We establish some general relations between Heegaard Floer based contact invariants. In particular, we observe that if the contact invariant of large negative, respectively positive, contact surgeries along a Legendrian knot does not vanish, then the Legendrian invariant, respectively the Legendrian inverse limit invariant, of that knot is non-zero. We use sutured Floer homology, and the limit constructions due to Golla, and Etnyre, Vela-Vick and Zarev.
\end{abstract}


\maketitle


\section{Introduction}
Heegaard Floer based contact invariants are the most used and powerful (though not complete) detector of tightness of contact manifolds. For Legendrian knots, they were independently developed from the point of view of sutured homology, knot Floer homology and grid homology, and how these invariants compare to each other attracted a lot of interest \cite{SV,BVV,G.i,EVVZ}, as did also their relationship to the invariants of contact surgeries along the knots \cite{Sah,LS.t,G.s,MT}. 

In contrast to the extensively studied behavior of surgeries on knots in the standard $3$-sphere, our interest is first in  surgeries along non-loose knots in overtwisted manifolds (that is, knots in overtwisted manifolds whose complement is tight). 

Recall that the contact $r$-surgery (whenever $r\neq\frac{1}{n}$) is not uniquely defined as it depends on choices of stabilizations of the Legendrian knot and its Legendrian push-offs (when described by surgery diagrams \cite{DG}). In the following, the special role will be played by the contact $r$-surgery with all stabilizations negative, and we will denote the resulting contact structure by $\xi_r^-$.

Based on the Legendrian surgeries, we can say the following about the Legendrian invariant of the knot and its torsion order.

\begin{thm}\label{thm:minus}
Let $L$ be a Legendrian knot in a contact manifold $(Y,\xi)$. If for every rational number $r\leq -1$ the contact $r$-surgery $\xi_r^-$ along $L$ has non-vanishing contact invariant, $c(Y_r(L),\xi_r^-)\neq 0$, then the Legendrian invariant of $L$ in $\HFK^-(-Y,L)$ is non-zero, $\Leg (L) \neq 0$.
\end{thm}

\begin{remark}
As we recall in Lemma \ref{lem:minus}, it suffices to check the non-vanishing condition for very negative integral surgeries.
\end{remark}

\begin{example}
Theorem \ref{thm:minus} allows us to prove the non-vanishing of the Legendrian invariants of non-loose $T_{(2,-2n+1)}$, conjectured by Lisca, Ozsv\'ath, Stipsicz and Szab\'o in \cite[Remark 6.11]{LOSSz}, as it is in greater generality carried out in \cite[Theorem 4.7]{M.n}.
\end{example}

Theorem \ref{thm:minus} is in fact interesting only when $c(Y,\xi)=0$, and then for integral surgeries we additionally observe the following.

\begin{prop}\label{prop:order}
Let $L$ be a Legendrian knot in a contact manifold $(Y,\xi)$ with $c(Y,\xi)=0$. If for all $n\geq m$ the contact $(-n)$-surgery along $L$ with $m$ positive stabilizations has non-zero contact invariant, then also the contact $(-n)$-surgery $\xi_{-n}^-$ has non-zero invariant; moreover, $m$ is at most the torsion order of $\Leg(L)$, that is $U^m\cdot\Leg(L)\neq 0$.
\end{prop}

\begin{example}
However, notice that taking any contact surgery in Theorem \ref{thm:minus} would not work. There are Legendrian knots for which contact $r$-surgery with all stabilizations positive results in a contact manifold with the non-zero invariant for all $r\leq -1$, but the Legendrian invariant of the knot vanishes. Such examples are the non-loose Legendrian negative torus knots $T_{(p,-q)}$ with $tb\leq-pq$ whose transverse approximation is loose (see \cite{M.n} for details).
\end{example}

\begin{example}
Inverse of Theorem \ref{thm:minus} is generally not true, not even when $\Leghat(L)\neq 0$. Take for an example a non-loose Legendrian right-handed trefoil $T_{(2,3)}$ in the overtwisted $(S^3,\xi)$ with Hopf invariant $d_3(\xi)=-1$ (see \cite[Figure 51]{EVVZ}), or more generally a non-loose Legendrian $T_{(2,2n+1)}$ in the overtwisted $(S^3,\xi)$ with Hopf invariant $d_3(\xi)=-2n+1$, that is, the knot $L(n)$ in \cite[Figure 9]{LOSSz}. As observed in \cite[Remark 6.5]{LOSSz}, some negative surgery on $L(n)$ produces a necessarily overtwisted contact structure on $S^3_{2n-1}(T_{(2,2n+1)})$, even though $\Leghat(L(n))\neq0$ according to \cite[Proposition 6.2]{LOSSz}.
\end{example}

In the case of positive surgeries, on the other hand, we provide an alternative view on and a generalization to all manifolds of some known results in the $3$-sphere, making use of the Legendrian inverse limit invariant. For the latter, in turn, we observe that it is not an independent invariant of Legendrian knots.

\begin{thm}\label{thm:plus}
Let $L$ be a Legendrian knot in a contact manifold $(Y,\xi)$. If there exists $R\geq 1$ such that for every $r>R$ the contact $r$-surgery $\xi_r^-$ along $L$ has non-vanishing contact invariant, $c(Y_r(L),\xi_r^-)\neq 0$, then the Legendrian inverse limit invariant of $L$ is non-zero, $\EHinverse (L)\neq 0$.
\end{thm}

\begin{remark}\label{rmk:1+}
As we write out in Lemma \ref{lem:plus}, it suffices to find a positive integral surgery with the non-vanishing invariant. In particular, it would suffice to choose the negative stabilizations only for the initial integral surgery. Note however that contact surgery with all stabilizations negative corresponds to inadmissible transverse surgery \cite{C}.
\end{remark}

\begin{prop}\label{prop:EHinverse}
For a Legendrian knot $L$ in a contact manifold $(Y,\xi)$, the non-vanishing of the inverse limit invariant $\EHinverse (L)\neq 0$ is equivalent to the non-vanishing of both the Legendrian hat invariant $\Leghat(L)\neq0$ and the ambient contact invariant $c(\xi)\neq0$.
\end{prop}

\begin{cor}\label{cor:tau_xi}
Let $L$ be a Legendrian knot in a contact manifold $(Y,\xi)$. If any positive (integral) contact surgery along $L$ results in a contact manifold with the non-vanishing contact invariant, then $(Y,\xi)$ is tight with $c(\xi)\neq 0$. When $L$ is null-homologous, both the invariants $\Leghat(L)$ and $\Leg(L)$ are non-zero, and $\Leg(L)$ generates one of the $\FF[U]$-towers of $\HFK^-(-Y,L)$; the Legendrian knot $L$ satisfies the Bennequin-type equality \[\tb(L)- \rot(L, [S]) = 2\tau_\xi(Y,L, [S])-1\] with respect to any Seifert surface $S$.
\end{cor}

\proof[Notation] Recall that the rank of $\HFK^-(-Y,\mathbf t_\xi, L)$ as $\FF[U]$-module is equal to the dimension of $\HFhat(-Y,\mathbf t_\xi)$, and that Hedden \cite[Definition 23]{H}  defines $\tau_\xi$ as the top grading of the tower corresponding to $c(\xi)$.

For example, in the case of the $3$-sphere, when looking at the knot in $-S^3$ corresponds to looking at the mirror knot in $S^3$, the invariant $\tau_\xi(S^3,L) = \tau(L) = -\tau(m(L))$.

\proof
Knowing Proposition \ref{prop:EHinverse}, both statements follow from Theorem \ref{thm:plus}.

The equality is obtained by the computation of the Alexander grading of $\Leghat(L)$ in terms of the classical invariants, as given by Ozsváth and Stipsicz in \cite[Theorem 1.6]{OS}.
\endproof

\begin{remark}
Alternatively, the first statement of Corollary \ref{cor:tau_xi} is obvious from the naturality of contact invariants with respect to positive surgeries, see \cite[Theorem 1.1]{MT} of Mark and Tosun. Meanwhile the second statement can be proven using surgery formulae, in the same way as \cite[Theorem 1.2]{MT} of Mark and Tosun for knots in $(S^3,\xi_\text{std})$.
\end{remark}

Recall that in \cite{LS.t} Lisca and Stipsicz defined an invariant of transverse and Legendrian knots $\tilde{c}$ using positive contact surgeries. It is defined as the class of the vector $(c(\xi_n^-(L)))_{n\in\NN}$ in the inverse system 
\[\left(\{\HFhat(-Y_n(L))\}_{n\in\NN}, \{\Phi_{n,m}=F_{\overline{W}_{\!n}}\circ\cdots\circ F_{\overline{W}_{\!m+1}}\}_{m<n}\right)\]
 of Heegaard Floer groups $\HFhat(-Y_n(L))$ of surgeries and the cobordism maps $F_{\overline{W}_{\!n}}: \HFhat(-Y_n(L)) \rightarrow \HFhat(-Y_{n-1}(L))$, defined through surgery exact triangles.

Theorem \ref{thm:plus} tells that non-vanishing of $\tilde{c}(L)$ implies non-vanishing of $\EHinverse (L)$. That answers \cite[Question 2]{EVVZ} of Etnyre, Vela-Vick and Zarev about the relationship between their inverse limit invariant $\EHinverse$ and the transverse invariant $\tilde{c}$ of Lisca and Stipsicz.

\begin{example}
Considering the inverse of Theorem \ref{thm:plus}, note that solely a non-zero $\Leghat(L)$ and the ambient contact manifold $(Y,\xi)$ having non-zero contact invariant, do not suffice for the non-vanishing of the contact invariant of large positive surgeries, as can be read from the conditions in \cite[Theorem 1.1]{G.s} of Golla  for the knots in $(S^3,\xi_\text{std})\ $(see also \cite[Theorem 1.2]{MT} of Mark and Tosun).

In particular, in the standard $3$-sphere $\EHinverse(L)\neq0$ is not equivalent to $\tilde{c}(L)\neq0$; for the latter, we additionally need that the underlying smooth knot satisfies the equality $\tau(L)=\nu(L)$.
\end{example}

\subsubsection*{Overview} The organization of the paper is straightforward. In Section \ref{Sec2} we briefly review relevant points about contact surgery and the Legendrian knot invarians defined in Heegaard Floer theory; throughout, we expect the basic knowledge of contact topology, in particular the convex surface theory, and the standard background in knot Floer and sutured Floer homology. In Section \ref{Sec3} we give the proofs of Theorem \ref{thm:minus}, Theorem \ref{thm:plus} and Proposition \ref{prop:EHinverse}.

\subsubsection*{Acknowledgement} I gained a lot from the work of Marco Golla and our conversation. I am grateful to Alberto Cavallo and Daniele Celoria for some clarifications in knot Floer homology. This work was supported by the European Research Council (ERC) under the European Union's Horizon 2020 research and innovation programme (grant agreement No 674978); it was carried out at the University of Oxford and I thank Andr\'as Juh\'asz for his support.

\section{Preliminaries}\label{Sec2}
\subsection{Contact surgery}
Recall that contact $r$-surgery (for any non-zero $r$) is performed along a Legendrian knot with the surgery coefficient $r$ measured relative to the contact framing; in addition to ordinary surgery, it prescribes for the contact structure to be preserved in the complement of the knot, while the extension to the glued-up torus needs to be tight. The possible contact structures on the solid torus filling are determined by the pair of two slopes: the boundary slope (that is, the slope of the dividing set on the boundary) equal to the initial contact framing $0$, and the meridional slope given by the surgery coefficient $r$. They are listed in the Honda's classification \cite{Ho.I} in terms of the shortest counter-clockwise path from $0$ to $r$ in the Farey graph (annotated by $(0,\infty,-1)$ in $(1,-1,\text{i})$ of $S^1$): the successive fractions along this path correspond to the successive slopes of basic slices glued to the knot complement. Explicitly, for a sequence $r_0=0, r_1,\dots, r_k=r$ we layer the glued-up torus from outside in into $k-1$ basic slices with boundary slopes $\{r_{i-1}, r_i\}$ for $i=1,\dots,k-1$, and fill it in with a solid torus of the boundary slope $r_{k-1}$ and meridional slope $r_k$; for each basic slice we have a choice of the sign, while the final solid torus admits a unique tight structure.

In particular case when $n\in\ZZ$, the sequence is comprised of $0,-1,\dots, n+1,n$ for $n<0$, and $0,\infty,n$ for $n>0$.
 
On the other hand, Ding and Geiges \cite{DG} gave an algorithm to convert contact $r$-surgery into a sequence of contact $(\pm1)$-surgeries, which encodes the convex decomposition of the glued-up torus in a form of the surgery diagram. Let us recall: If the surgery coefficient is $r=\frac{p}{q}$ and $m\in\NN$ is the minimal such that $\frac{p}{q-mp}<0$, and $\frac{p}{q-mp}=[a_0,\dots,a_k]$, the slicing can be described on Legendrian push-offs of the surgered knot $L$ as follows. First we perform contact $(+1)$-surgery along $m$ push-offs of $L$, then for each successive $i^\text{th}$ continued fraction block we do $(-1)$-surgery along $L_i$ where $L_i$ is obtained from $L_{i-1}$ by Legendrian push-off and additional $a_i-1$ stabilizations, and $L_0=L$ stabilized $a_0-1$ times. All possible contact structures on glued-up torus are then covered by all possible choices of positive or negative stabilizations.
 
In the case $n\in\ZZ$, this amounts to a single $(-1)$-surgery along $(n-1)$-times stabilized Legendrian knot for $n<0$, and to a single $(+1)$-surgery along the knot followed by $(-1)$-surgery along $n-1$ of its once stabilized push-offs for $n>0$.
 
\subsection{Legendrian invariants in Heegaard Floer theory} The initial invariant of Legendrian knots was defined by Honda, Kazez and Mati\'c \cite{HKM.eh} in sutured Floer homology \cite{J}, in analogy to the contact invariant of closed manifolds due to Ozsváth and Szabó \cite{OSz.c}. Subsequently, Stipsicz and Vértesi \cite{SV}, Golla \cite{G.i}, and Etnyre, Vela-Vick and Zarev \cite{EVVZ} found interpretations of various knot Floer invariants in terms of the sutured Floer homology, which we aim to recall in this subsection.

\subsubsection*{Convention} In accordance with the notation in the previous subsection on contact surgery, we parametrize the boundary of the Legendrian knot complement by the meridian of the knot taking the slope $\infty$, and the contact framing taking the slope $0$. Thus, in the case of a nullhomologous knot, our slope $0$ corresponds to the slope $\tb$ with respect to the Seifert framing.

\subsubsection*{$\EH$-invariants and gluing maps} Honda, Kazez and Mati\'c \cite{HKM.eh} define $\EH(Y,\xi)$ of a contact manifold with convex boundary as a class in $\SFH(-Y,-\Gamma_\xi)$ where $\Gamma_\xi$ consists of dividing curves of $\xi$ on $\partial Y$. The contact invariant $c(Y,\xi)$ is then identified with $\EH((Y,\xi)\backslash(B^3,\xi_\text{std}))$ in $\SFH(Y(1))$, and in case of a Legendrian knot $L$, the invariant $\EH(L)$ is set to be $\EH(Y(L),\xi_L)$ in $\SFH(-(Y\backslash\nu L), -\Gamma_\text{0})$ where $\Gamma_\text{0}$ is a pair of oppositely oriented closed curves of slope $0$ on the boundary torus. The crucial property of these sutured contact invariants is their behaviour under gluing diffeomorphisms \cite{HKM.glue}: for sutured manifolds $(Y,\Gamma)\subset(Y',\Gamma')$, a contact structure $\zeta$ on $Y'\backslash Y$, compatible with $\Gamma\cup\Gamma'$, induces a map 
\[\Phi_\zeta: \SFH(-Y,-\Gamma) \rightarrow \SFH(-Y',-\Gamma')\]
which in case of contact manifolds with convex boundary connects the $\EH$-classes, that is
\[\Phi_\zeta(\EH(Y,\xi_Y))=\EH(Y',\xi_Y\cup\zeta).\]

Studying Legendrian knots, the sutured manifolds, we are specifically interested in, are the knot complements with various boundary slopes $(Y(L),\Gamma_s)$ and their (punctured) Dehn fillings $(Y_r(L))(1)$. The key gluing maps are:
\begin{itemize}[leftmargin=.5cm] 
\item $\phi_{s_0,s_1}: \SFH(-Y(L),-\Gamma_{s_0})\rightarrow \SFH(-Y(L),-\Gamma_{s_1})$, associated to the addition of a tight toric annulus $(T^2\times I,\zeta_{s_0,s_1})$; in particular, the maps $\sigma^\pm_{s_0,s_1}$ associated to the basic slices, and 
\item $\psi_{r(s)}: \SFH(-Y(L),-\Gamma_s)\rightarrow \HFhat(-Y_r(L))$, associated to the surgery $2$-handle, glued with slope $r$ to the torus of the boundary slope $s$; note that in case $s=0$ this corresponds to gluing-up a (punctured) tight solid torus $(V_r(1),\zeta_r)$, and  when $s$ and $r$ are connected in the Farey graph, in which case we denote it simply by $\psi_r$, the contact filling is unique.
\end{itemize}

\subsubsection*{Invariant $\Leghat$} Stipsicz and V\'ertesi \cite{SV} interpret the Legendrian invariant $\Leghat(L)\in\HFKhat(-Y,L)$, originally defined by Lisca, Ozsv\'ath, Stipsicz and Szab\'o \cite{LOSSz}, as the $\EH$-invariant of the complement of a Legendrian knot $(Y(L),\xi_L,\Gamma_\text{0})$ completed by the negative basic slice with boundary slopes $0$ and $\infty$. Therefore, if we denote the completion $\xi_L\cup\zeta^-_{0,\infty}$ by $\overline{\xi_\infty}$, we have
\[\Leghat(L)=\EH(Y(L),\overline{\xi_\infty}) \in\SFH(-Y(L),-\Gamma_\infty).\]

\subsubsection*{Invariant $\Leg$} Golla \cite{G.i}, explicitly for knots in the $3$-sphere,  and later Etnyre, Vela-Vick and Zarev \cite{EVVZ} give corresponding reinterpretation for the Legendrian invariant $\Leg(L)\in\HFK^-(-Y,L)$ \cite{LOSSz}. They both observe that $\SFH(-Y(L),\Gamma_{-i})$ together with the negative bypass attachments 
\[\phi^-_{-i,-j}=\sigma^-_{-i,-i-1}\circ\cdots\circ \sigma^-_{-j+1,-j}: \SFH(-Y(L),-\Gamma_{-i}) \rightarrow \SFH(-Y(L),-\Gamma_{-j})\] 
form a direct system $(\{\SFH(-Y(L),\Gamma_{-i})\}_{i\in\NN}, \{\phi^-_{-i,-j}\}_{j>i})$ and that its direct limit $\underrightarrow{\SFH}(-Y,L)$  is isomorphic to $\HFK^-(-Y,L)$. Furthermore, they observe that $\EH$-invariants of the negative stabilizations of $L$ respect this direct system, and that the class of the vector $(\EH(L^{i-}))_{i\in\NN}=\EHdirect(L)$ is taken to $\Leg(L)$ under the above isomorphism. 

\subsubsection*{Invariant $\EHinverse$} Etnyre, Vela-Vick and Zarev \cite{EVVZ} consider also a parallel inverse system with positive boundary slopes $(\{\SFH(-Y(L),\Gamma_{i})\}_{i\in\NN}, \{\phi^-_{j,i}\}_{j>i})$, whose inverse limit $\underleftarrow{\SFH}(-Y,L)$  is isomorphic to $\HFK^+(-Y,L)$. They define the inverse limit invariant $\EHinverse(L)$ to be the class of the vector $(\EH(Y(L),\overline{\xi_i}))_{i\in\NN}$ where $\overline{\xi_i}$ equals the extension by two negative basic slices $\xi_L\cup\zeta^-_{0,\infty}\cup\zeta^-_{\infty,i}$.

\section{Proofs}\label{Sec3}

\begin{lem}\label{lem:minus}
If the contact invariant $c(Y_r(L),\xi_r^-)$ is non-zero for all large negative integers $r$, then it is non-zero for all $r\leq -1$.
\end{lem}
\proof
First, utilizing contact surgery presentation of Ding and Geiges \cite{DG}, we know that $(Y_r(L),\xi_r^-)$ for $r\in (-n,-n+1)$ and $n\in \NN$ can be obtained by Legendrian surgery on $(Y_{-n}(L),\xi_{-n}^-)$; hence, if the integral surgeries have non-zero contact invariant, so do the rational ones. 

Furthermore, if the contact invariant of $(Y_{-n}(L),\xi_{-n}^-)$ is non-zero, the contact invariant of $(Y_{-m}(L),\xi_{-m}^-)$ for $m<n$ is non-zero too. Indeed, the contact $(-n)$-surgery along a Legendrian knot $L$ equals some contact $(-n+1)$-surgery along a once stabilized knot $L'$. And, since every stabilization of the knot lies on a (boundary-parallel) stabilization of the page compatible with the original Legendrian knot, the capping off morphism of Baldwin implies, that the contact invariant of the Legendrian surgery along a stabilized knot $L'$ vanishes once the contact invariant of the Legendrian surgery along the knot $L$ is zero  \cite[Theorem 1.7]{B}. 
\endproof 

\proof[Proof of Theorem \ref{thm:minus}]
Under the identification of $\Leg(L)$ with $\EHdirect(L)$, we need to show that $\EH(L^{i-})$ remains non-zero as the number of negative stabilizations goes to infinity.

Of course, the complement $(Y(L),\xi_L)$ of the standard neighborhood of the Legendrian knot $L$, as well as the complements \[(Y(L^{i-}),\xi_{L^{i-}}) = (Y(L),\xi_L\cup\bigcup_{j=-1}^{-i}\zeta^-_{j+1,j})\] of its stabilizations embed into $(Y,\xi)$ as sutured submanifolds. Hence, the induced map on the sutured Floer homology 
\[\psi_{\infty(-i)}: \SFH(-Y(L), -\Gamma_\text{-i}) \rightarrow \SFH(-Y(1))\] 
takes
\[\psi_{\infty(-i)} (\EH(L^{i-}))=c(Y,\xi).\]

Negative contact surgery takes away some of the twisting around the knot; indeed, we can always replace it by an admissible transverse surgery, which in turn can be understood as a particular contact cut (see \cite{BE} for details). However, in contact $(-n)$-surgery $(Y_{-n}(L),\xi_{-n}^-)$ still embed the standard complements of $L^{i-}$ for all $i<n$, or more precisely, contact $(-n)$-surgery $(Y_{-n}(L),\xi_{-n}^-)$ equals the complement of $L^{(n-1)-}$ filled by the solid torus of the boundary slope $-n+1$ and meridional slope $-n$. Therefore, we have
\[\psi_{-n(-i)}:\SFH(-Y(L),-\Gamma_{-i}) \rightarrow \SFH((-Y_{-n}(L))(1))\]
and
\[\psi_{-n(-i)}(\EH(L^{i-}))=c(Y_{-n}(L),\xi_{-n}^-) \text{ for all } i<n.\]

Since we assumed that the contact invariants of large negative contact surgeries with all stabilizations negative do not vanish, the $\EH$-invariants of the negative stabilizations of the Legendrian knot $L$ do not vanish either, and thus $\EHdirect(L)=\Leg(L)\neq 0$.
\endproof

\proof[Proof of Proposition \ref{prop:order}]
The condition that the contact $(-n)$-surgery (for every $n\geq m$) along $L$ with $m$ positive stabilizations has non-zero contact invariant, means that there is a contact $-(n-m)$-surgery $\xi_{m-n}^-$ along $L^{m+}$ with non-zero contact invariant. Then, as a consequence of the capping off morphism of Baldwin \cite[Theorem 1.7]{B}, there is also a contact $-(n-m)$-surgery $\xi_{m-n}^-$ along $L$ with non-zero contact invariant for every $n\geq m$. 

On the other hand, when $c(Y,\xi)=0$, we know that the $U$-order of $\Leg(L)$ is finite, say $o$; in particular, the invariant of the $o$-times positively stabilized $L$ vanishes. Now, since for all $n$ a contact $-(n-m)$-surgery along $L^{m+}$ has non-zero contact invariant,  the invariant $\Leg(L^{m+}) =U^m\cdot\Leg(L) \neq0$ by Theorem \ref{thm:minus}, and hence $m<o$.
\endproof

\begin{lem}\label{lem:plus}
If there is $n\in\NN$ such that $c(Y_n(L),\xi_n^-)\neq0$, then the contact invariant $c(Y_r(L),\xi_r^-)\neq0$ for all $r\geq n$.
\end{lem}
\proof
According to the contact surgery presentation of Ding and Geiges \cite{DG}, every $(Y_r(L),\xi_r^-)$ for $r\geq n$ and $n\in \NN$ can be obtained by Legendrian surgery on $(Y_{n}(L),\xi_n^-)$. 
\endproof

\proof[Proof of Theorem \ref{thm:plus}]
Opposed to the complements of the standard neighborhoods of the Legendrian knot stabilizations $L^{i-}$, the completion of the knot complement $(Y(L),\overline{\xi_\infty})$ never embeds in $(Y,\xi)$, and hence neither do its extensions $(Y(L),\overline{\xi_i})$ for $i\in\NN$. 

Nevertheless, they do embed in positive contact surgeries with all stabilizations negative; indeed, in contrast to the negative surgeries these correspond to inadmissible transverse surgery and can be thought of as adding twisting along a Legendrian knot (see \cite{C}). Explicitly, in contact $r$-surgery where $r>0$, embed all the integral extensions for $n>r$. In particular, for $r=\frac{2n-1}{2}$ the contact $r$-surgery equals the $n$-extension $(Y(L),\overline{\xi_n})$ filled by the solid torus of the boundary slope $n$ and meridional slope $r$, and the morphism induced by the surgery $2$-handle
\[\psi_r:\SFH(-Y(L),-\Gamma_n) \rightarrow \SFH((-Y_r(L))(1))\]
takes
\[\psi_r(\EH(Y(L),\overline{\xi_n})) = c(Y_r(L),\xi_r^-).\]

Now, as we have assumed to have at least one positive surgery with non-vanishing contact invariant, all larger surgeries -- as observed in Remark \ref{rmk:1+} -- have non-vanishing invariant as well. Hence, for big enough $n$ the $n$-extension embeds into closed contact manifold with non-zero contact invariant, and has non-zero $\EH$. Thus, $\EHinverse(L)\neq 0$.
\endproof

\proof[Proof of Proposition \ref{prop:EHinverse}]
According to Etnyre, Vela-Vick and Zarev \cite[Theorem 1.7]{EVVZ}, there is a commutative triangle
\[
\begin{tikzpicture}
\matrix (m) [matrix of math nodes, row sep=3em,column sep=3em]
{ \underleftarrow{\SFH}(-Y,L) &    & \HFK^+(-Y,L) \\
              &\ \  \HFKhat(-Y,L)  &             \\ };
\path[->,font=\small]
	(m-1-1)	edge node[above] {$I_+$} (m-1-3);
\path[<-,font=\small]
	(m-1-1) edge node[auto,swap] {$\Phi_{\text{dSV}}$} (m-2-2)
	(m-1-3) edge node[auto] {$\iota_*$} (m-2-2);
\end{tikzpicture}
\]
where $I_+$ is the isomorphism and $\iota_*$ is the map induced on homology by the inclusion of complexes. Additionally, by \cite[Theorem 1.7]{EVVZ}, the map $\Phi_\text{dSV}$ sends the Legendrian invariant $\Leghat(L)$ into $\EHinverse(L)$.

Now, recall that $\HFK^+$ is defined as the homology of the graded object associated to $\CFK^\infty$ quotient by $U\cdot\CFK^-$:
\[ \CFK^+:=\CFK^\infty/\ U\cdot\CFK^- \text{ and } \HFK^+:=\text{H}_*(\gCFK^+).  \]
Hence, the map $\iota_*$ sends to zero exactly those elements of $\HFKhat$ which correspond to torsion elements in $\HFK^-$.

Existence of the map $\Phi_\text{dSV}$ requires that $\Leghat(L)\neq 0$ when $\EHinverse(L)\neq 0$. Moreover, since the Legendrian invariant $\Leghat(L)$ as an element of $\HFK^-(-Y,L)$ is non-torsion if and only if the ambient contact structure has non-zero invariant \cite[Theorem 1.2]{LOSSz}, the result follows.
\endproof


\end{document}